\documentclass[12pt]{article}
\usepackage{amsfonts, amsmath, amsthm}

\newtheorem{theorem}{Theorem}[section]
\newtheorem{lemma}[theorem]{Lemma}

\newtheorem{cor}[theorem]{Corollary}
\newtheorem{prop}[theorem]{Proposition}

\def\QQ{\mathbb{Q}}
\def\RR{\mathbb{R}}

\def\calO{\mathcal{O}}

\def\alg{\mathrm{alg}}
\def\an{\mathrm{an}}
\def\con{\mathrm{con}}

\def\imm{\mathrm{imm}}
\def\sep{\mathrm{sep}}
\def\perf{\mathrm{perf}}

\def\beq{\begin{equation}}
\def\eeq{\end{equation}}

\def\fp{\frac{1}{p}}

\def\Galg{\Gamma^{\alg}}
\def\Galgcon{\Galg_{\con}}
\def\Galgancon{\Galg_{\an,\con}}

\def\Gancon{\Gamma_{\an,\con}}
\def\Gcon{\Gamma_{\con}}
\def\Gimm{\Gamma^{\imm}}
\def\Gimmcon{\Gimm_{\con}}
\def\Gimmancon{\Gimm_{\an,\con}}

\def\Gsep{\Gamma^{\sep}}
\def\Gsepcon{\Gsep_{\con}}

\def\Oalg{\Omega^{\alg}}
\def\Oan{\Omega_{\an}}
\def\Oimm{\Omega^{\imm}}

\def\bv{\mathbf{v}}
\def\bw{\mathbf{w}}
\def\bx{\mathbf{x}}
\def\by{\mathbf{y}}

\begin{document}

\title{The Newton polygons of overconvergent $F$-crystals}
\author{Kiran S. Kedlaya}
\date{May 31, 2001}

\maketitle

\begin{abstract}
It is conjectured that every overconvergent $(F, \nabla)$-crystal
over $k((t))$ is potentially semistable (equivalently, quasi-unipotent),
and so has ``generic''
and ``special'' Newton polygons. It is easy to construct a Newton polygon
for an arbitrary overconvergent $F$-crystal that coincides with the generic
Newton polygon for potentially semistable crystals. We give an analogous
construction for the special Newton polygon. 
In a subsequent preprint, we use this construction to prove the
aforementioned conjecture.
\end{abstract}

Crew \cite{bib:crew1} asked whether every overconvergent $(F, \nabla)$-crystal
over $k((t))$ is quasi-unipotent.
An equivalent conjecture \cite[Conjecture~4.12]{bib:me3}
is that
every such crystal is potentially semistable,
which is to say, isomorphic to a log-crystal over $k[[u]]$ for some finite
extension $k((u))$ of $k((t))$.
(A related conjecture is the ``$p$-adic monodromy'' conjecture of
Fontaine \cite{bib:fon}, that every de~Rham representation is potentially
semistable.)
%A related conjecture is the
%``monodromy conjecture'' of Berger and Colmez \fixme{get reference}.
Under this conjecture, every 
overconvergent $(F, \nabla)$-crystal over $k((t))$
has two well-defined Newton 
polygons, corresponding to the fibres of the crystal over the generic
and special points of $k[[u]]$. By Grothendieck's specialization theorem,
these polygons have the same endpoints and the special polygon lies on or
above the generic polygon.

One difficulty in proving this conjecture has been the lack of an \emph{a
priori} description of the special Newton polygon. The generic Newton polygon
coincides with the Newton polygon as a crystal over $\overline{k((t))}$.
The special Newton polygon is only uniquely determined from this data alone
when the generic slopes
are all equal, in which case the generic and special Newton polygons must
coincide. In fact, the conjecture was previously
known to hold when the generic slopes
are all equal, by a result of Tsuzuki \cite{bib:tsu1}.

In this paper, we describe a construction of a Newton polygon for an arbitrary
overconvergent $F$-crystal over $k((t))$ which coincides with the special
Newton polygon for potentially semistable crystals and has some of the
expected formal properties of the special Newton polygon. In a subsequent
preprint \cite{bib:me6}, we use this construction to prove Crew's
conjecture. (Note that Yves Andr\'e and Zoghman Mebkhout have each 
independently announced proofs of this conjecture, at least in the case
where the residue field is the algebraic closure of a finite field.)
%This allows us
%to extend some structural results about $F$-crystals over $k[[t]]$ to
%the overconvergent setting, including Katz's decomposition theorem.

\section{Notations}
We maintain the notations of \cite{bib:me3}, recalled below.
(Note:
wherever it appears in the table, $*$ represents an unspecified
decoration.)

\begin{list}%
{}{\setlength{\itemsep}{0pt}
\setlength{\parsep}{0pt}
\setlength{\labelwidth}{0.75in}
\setlength{\leftmargin}{1in}
\setlength{\labelsep}{0.25in}}
\item[$k$]
An algebraically closed field of characteristic $p>0$
\item[$\calO$]
A finite extension of the Witt ring $W(k)$
\item[$\sigma$ on $\calO$]
An automorphism of $\calO$ lifting the absolute Frobenius on $k$
\item[$\calO_0$]
The elements of $\calO$ fixed by $\sigma$
\item[$|\cdot|$]
The valuation of $\calO$ normalized so that $|p| = p^{-1}$
\item[$K$]
The field of formal power series over $k$
\item[$K^{\perf}$]
The perfect closure of $K$
\item[$K^{\sep}$]
The separable closure of $K$
\item[$K^{\alg}$]
The algebraic closure of $K$
\item[$K^{\imm}$]
The ring of series $x = \sum_{i \in I} x_i t^i$ over $k$ with $I \subseteq \QQ$ well-ordered
\item[$\Omega_t$ (= $\Omega$)]
The power series ring $\calO[[t]]$
\item[$\Gamma$]
The $p$-adic completion of $\Omega[t^{-1}]$
\item[$\sigma$ on $\Gamma$]
An endomorphism lifting $x \mapsto x^p$ compatible 
with $\sigma$ on $\calO$
\item[$\sigma_t$]
The endomorphism of $\Gamma$ with $t \mapsto t^p$
\item[$\Gamma^L$]
The $p$-adically complete
extension of $\Gamma$ with residue field $L$
\item[$\Gamma^*$]
Equal to $\Gamma^{K^*}$ for $* \in \{\perf, \sep,
\alg, \imm\}$
\item[$\Gamma^{\alg(c)}$]
The subring of $x = \sum_i x_i t^i \in \Galg$
for which for each $n \geq 0$, there exists $r_n$
such that $|x^{\sigma^{n}}-r_n|
< p^{-cn}$
\item[$\Gamma^{*}_{\con}$]
The ring of $x = \sum_{i=-\infty}^{\infty} x_i t^i \in \Gamma^*$ with
$x_i \in \calO$ and  $\liminf_{i \to \infty} v_p(x_{-i})/i > 0$
\item[$\Gancon$]
The ring of $x = \sum_{i=-\infty}^{\infty} x_i t^i$ with
$x_i \in \calO[\fp]$, $\liminf_{i \to \infty} v_p(x_{-i})/i > 0$,
and $\liminf_{i \to +\infty} v_p(x_i)/i \geq 0$ 
\end{list}

We also mention two rings that are not explicitly defined in 
\cite{bib:me3}. The ring $\Oimm$ consists of those elements $x = 
\sum x_i t^i$ of $\Gimm$ with $x_i=0$ for $i<0$. The ring $\Oalg$
is the intersection of $\Oimm$ with $\Galg$.

\section{Analytic rings revisited}
In this section, we make a more careful study of some of the 
rings introduced in \cite{bib:me3} but not used extensively there.
In particular, we need a careful analysis of the rings $\Gancon$ and its
extensions.

Recall that the ring $\Gancon$ is defined as the set of series
$x = \sum_{n=-\infty}^{\infty} x_n t^n$ with $x_n \in \calO[\fp]$
satisfying $v_p(x_n) \geq -cn$ for $n$ sufficiently negative (and some
choice of $c$) and $v_p(x_n) = o(n)$ for $n$ sufficiently positive.
Also recall that the ring $\Gimm$ is defined as the set of series
$x = \sum_{n \in \QQ} x_n t^n$ with $x_n \in \calO$ such that for each
$r>0$, the set of $n$ with $|x_n| > p^{-r}$ is a well-ordered subset of
$\QQ$. To put these together, we define the ring $\Gimmancon$ as the
set of series $x = \sum_{n \in \QQ} x_n t^n$, with $x_n \in \calO[\fp]$,
satisfying the following
conditions:
\begin{enumerate}
\item
For each $r \in \RR$, the set of $n \in \QQ$ with $|x_n| > p^r$ is well-ordered.
\item 
There exists a constant $c>0$ such that for $n$ sufficiently negative,
$|x_n| < p^{cn}$.
\item
For every constant $d>0$, we have $|x_n| > p^{dn}$ for $n$ sufficiently
positive.
\end{enumerate}

This ring has an unusual property not shared by any of the rings introduced
so far: it contains nonzero solutions of the equation $x^\sigma = \lambda x$
for $\lambda \in \calO$ not a unit. For example, if $\lambda^\sigma
= \lambda$, then $x = \sum_{n=-\infty}^\infty
\lambda^{-n} t^{p^n}$ is a solution.

Notice that in $\Gimmancon$, the equation $x^\sigma = \mu x$ with
$\mu \in \calO$ has nontrivial solutions whenever $|\mu| \leq 1$. (By
contrast, in $\Gimmcon$, there are no solutions unless $|\mu| = 1$.) In
fact, it is easy to write down all such $x$: they are given as
\[
x = \sum_{n=-\infty}^\infty \mu^n y^{\sigma^{-n}}, \qquad
y = \sum_{i \in [a, ap)} y_i t^i,
\] 
where $a$ is a fixed positive rational.

Recall that we view $\Galgcon$ as a subring of $\Gimmcon$.
Within $\Gimmancon$, we can correspondingly identify the subring $\Galgancon$
consisting of elements $x = \sum x_i t^i$ such that for each $j$,
$\sum_{i<j} x_i t^i \in \Galgcon[\fp]$. Alternatively, for each $j$,
there exists $y \in \Galgcon[\fp]$ with $x_i = y_i$ for $j<i$;
the equivalence of these two conditions reduces to the fact that the
truncation of an element of $\Galgcon$ is still in $\Galgcon$.
This statement can be proved directly or
deduced from the classification of algebraic generalized
power series \cite[Theorem~8]{bib:me1}.

\begin{lemma} \label{lem:rankone}
Suppose $m$ is a positive integer,
$\lambda,\mu \in \calO$ are nonzero 
and $x = \sum x_i t^i \in R$, for $R = \Galgancon$ or $R=\Gimmancon$.
\begin{enumerate}
\item[(a)]
If $|\lambda| \geq |\mu|$, then there exists $y \in R$
such that $\lambda y^{\sigma^m} - \mu y = x$.
\item[(b)]
If $|\lambda| \leq |\mu|$ and $x_i=0$ for $i<0$, then there exists
$y \in R$ such that $\lambda y^{\sigma^m} - \mu y = x$.
\item[(c)]
If $x \in \Gimmcon$, $|\lambda| \leq |\mu|$ and there exists
$y \in R$ such that $\lambda y^{\sigma^m} - \mu y = x$, then
$\mu y \in \Gimmcon$.
\end{enumerate}
\end{lemma}
\begin{proof}
There is no loss of generality in assuming $\lambda, \mu \in \calO_0$.

(a)
First suppose $|\lambda| = |\mu|$; without loss of generality, we 
may assume $\lambda = \mu = 1$.
Write $x = a+x_0+b$ with $a = \sum_{i<0} x_i t^i$ and $b = \sum_{i > 0}
x_i t^i$.
Set
\[
c = \sum_{i<0} \sum_{n=1}^\infty x_{ip^{mn}}^{\sigma^{-mn}}, \qquad
d = \sum_{i>0} \sum_{n=0}^\infty x_{ip^{-mn}}^{\sigma^{mn}}.
\]
The first inner sum converges because $v_p(x_i) \geq O(-i)$ for $i \to 
-\infty$, while the second inner sum converges because $|x_i| \to 0$ as $i$
runs over any decreasing sequence. It is easily checked that
$\lambda c^{\sigma^m} - \mu c = b$ and that
$\lambda d^{\sigma^m} - \mu d = b$. Let $y_0 \in \calO$ be a solution
of $y_0^\sigma - y_0 = x_0$; then we may set $y = c+d+y_0$ to obtain
a solution of
$\lambda y^{\sigma^m} - \mu y = x$.

Next, suppose $|\lambda| > |\mu|$.
Write $x = a+b$ with $a = \sum_{i<1} x_i t^i$ and $b = \sum_{i \geq 1}
x_i t^i$. (The index 1 could be replaced by any positive index.)
We wish to set
\[
c = \sum_{n=1}^\infty \frac{a^{\sigma^{-mn}} \mu^{n-1}}{\lambda^n}, \qquad
d = \sum_{n=0}^\infty -\frac{b^{\sigma^{mn}} \lambda^n}{\mu^{n+1}}
\]
The first sum converges $p$-adically in $\Galgcon[\fp]$ or
$\Gimmcon[\fp]$ to a solution of $\lambda c^{\sigma^m} - \mu c = a$.
The second sum
converges $t$-adically (since $b^{\sigma^{mn}}$ has no coefficients
of index less than $p^{mn}$) to a solution of
$\lambda d^{\sigma^m} - \mu d = b$.
Thus we may set $y=c+d$ to obtain a solution of
$\lambda y^{\sigma^m} - \mu y = x$.

(b)
Set
\[
y = -\sum_{i\geq 0} \sum_{n=0}^\infty x_{ip^{-mn}}^{\sigma^{mn}}
\frac{\lambda^{n}}{\mu^{n+1}};
\]
the inner sum converges because $x_{ip^{-mn}}$ is bounded and 
$(\lambda/\mu)^n \to 0$.
Then it is easily verified that $\lambda y^{\sigma^m} - \mu y = x$.

(c)
Suppose $y \notin \Gimmcon$; put $y = \sum_i y_i t^i$ and let $j$ be
the smallest index such that $|\mu y_j| > 1$. First suppose $j \geq 0$.
Comparing the coefficients
of $t^j$ in the equation $\lambda y^{\sigma^m} - \mu y = x$, we
have $\lambda y_{j/p^m}^{\sigma^m} - \mu y_j = x_j$. In this equation,
$|\mu y_j| > 1$ but $|x_j| \leq 1$, so $|\lambda y_{j/p^m}^{\sigma^m}| 
= |\mu y_j| > 1$.
On the other hand, $|\lambda y_{j/p^m}^{\sigma^m}| \leq
|\mu y_{j/p^m}| \leq 1$ since $j/p^m < j$, contradiction.

Now suppose $j < 0$. In this case, the above argument gives
$|y_{j/p^m}| = |y_j \mu/\lambda|$. By comparing the coefficients of
$t^{j/p^m}, t^{j/p^{2m}}, \dots$ as well, we obtain
$|y_{j/p^{lm}}| = |y_j \mu^l/\lambda^l|$ by induction on $l$. But this
conclusion contradicts the fact that $|y_i|$ is bounded for $i < 0$.
\end{proof}

\section{Factorization theorems}

In this section, we pick up the thread begun in \cite[Section~4]{bib:me3},
and prove some additional
factorization theorems for analytic rings. Our goal is to prove
that finitely generated ideals in analytic rings are principal
(Lemma~\ref{lem:principal}),
a fact which will be crucial in the proof of the main theorem.

We begin with a lemma that will allow us to focus on $\Oan$ and its
analogues instead of on $\Gancon$ and its analogues. In fact, the
corresponding statement for $\Gancon$ and $\Oan$ is \cite[Lemma~4.7]{bib:me3},
and the proof here is simply a careful generalization of the proof there.
\begin{lemma} \label{lem:immfact1}
Every nonzero element of $\Gimmancon$ (resp.\ $\Galgancon$)
can be factored as the product
of an element of $\Gimmcon$ and an element of $\Oimm_{\an}$
(resp.\ an element of $\Galgcon$ and an element of $\Oalg_{\an}$),
the latter having nonzero constant coefficient.
\end{lemma}
\begin{proof}
Let $x = \sum_i x_i t^i$
be an element of $\Gimmancon$. (The proof for $\Galgancon$ is
analogous, so we omit reference to it hereafter.) Let $c$ be an irrational
number, so that $\min_i \{v_p(x_i) + ci\}$ occurs for a unique value of $i$.
Without loss of generality, we may assume that value is $i=0$ and that
$x_0 = 1$. (Otherwise,
we can multiply by a suitable constant times a suitable power of $t$.)
Put
\[
r = \min_{i<0} \left\{ \frac{v_p(x_i)}{-i} \right\}, \qquad
s = \min_{i>0} \left\{ \frac{-v_p(x_i)}{i} \right\};
\]
then by construction, $r<c < s$.

Now define a sequence $\{y_n\}_{n=0}^\infty$ as follows.
Let $k$ be the smallest index such that $v_p(x_k) < 0$. (If no such 
$k$ exists, then $x \in \Gimmcon$ and there is nothing to prove.) Set $y_0=0$.
To define $y_{n+1}$ from $y_n$, set $(1+y_n) x = \sum_i y_{n,i} t^i$ and let
\[
y_{n+1} = -1 + (1+y_n) \left( \sum_{i < k} y_{n,i}t^i \right)^{-1}.
\]
We will show that $\{y_n\}$
converges,
 in a suitable sense, to a limit $y$ and that $yx \in \Oimm_{\an}$.

Suppose $d \geq  0$ satisfies $v_p(y_{n,i}) \geq -ri + d$ for all $i<k$.
Let $a_n = -1 + \sum_{i < k} y_{n,i} t^i$,
$b_n = \sum_{i \geq k} y_{n,i} t^i$ and $c_n = (1+a_n)^{-1} - 1$.
Then $c_n = \sum_{j=1}^\infty (-a_n)^j$; if we put $c_n = 
\sum_i c_{n,i} t^i$, it follows that $v_p(c_{n,i}) \geq -ri + d$ for all
$i<k$. Then
\begin{align*}
(1+y_{n+1})x &= (1+y_n)x(1+b_n) \\
&= (1 + a_n + b_n) (1 +c_n) \\
&= 1 + b_n(1 + c_n).
\end{align*}
Therefore for $j < k$ nonzero,
\begin{align*}
v_p(y_{n+1,j}) &= v_p\left(\sum_{i\geq k} y_{n,i} c_{n,j-i} \right) \\
&\geq \min_{i \geq k} \{ v_p(y_{n,i}) + v_p(c_{n,j-i}) \} \\
&\geq \min_{i \geq k} \{ -si - r(j-i) + d \} \\
&= \min_{i \geq k} \{ -rj + d + (r-s)j \} \\
&= -rj + d + (r-s)k.
\end{align*}
For $n=0$, the initial inequality holds for $d=0$. Therefore for $j<k$,
$v_p(y_{n,j}) \geq -rj + n(r-s)k$ by induction on $n$, and the
same inequality holds for the $c_{n,j}$ as noted above.

Let $R$ be the subring of $\Gimmcon$ consisting of series $z = \sum_i z_i t^i$
such that $v_p(z_i) \geq -ri$ for all $i<0$. Then $R$ carries a valuation
$v'$ defined by $v'(z) = \min_i \{ v_p(z_i) + ri \}$. With respect to $v'$,
the product $\prod_{n =1}^\infty (1 + c_n)$ converges; we call
the limit $y$. Moreover, we can extend $v'$ to the subring of 
$\Gimmancon$ defined by $v_p(z_i) \geq -ri$ for all $i$,
in which it is clear that $y_n x \to yx$. In particular, for $j<k$,
$(yx)_j = \lim_{n \to \infty} y_{n,j} = 0$.
Thus $yx \in \Oimm_{\an}$, and we may factor $x$ as $y^{-1}(yx)$ with
$y^{-1} \in \Gimmcon$ and $yx \in \Oimm_{\an}$, as desired.
\end{proof}

Unless otherwise specified, throughout this section $R$ will denote
one of $\Oan$, $\Oalg_{\an}$, or $\Oimm_{\an}$.
For each $r > 0$, define the norm $w_r$ as follows:
\[
w_r(x) = \max_{i \geq 0} \{ |x_i| p^{-ri} \}
\qquad x\neq 0.
\]
We can use these norms to endow $R$ with a Fr\'echet topology. Recall that
this means a sequence $\{x_n\}$ is Cauchy if and only if for
each $r>0$ and $\epsilon > 0$,
there exists $N$ such that $w_r(x_m - x_n)<\epsilon$ for $m,n \geq N$.
Note that each $w_r$ extends to
a certain subring of $\Gimmancon$, but not to the whole ring.

\begin{prop}
The ring $R$ is complete for the Fr\'echet topology.
\end{prop}
\begin{proof}
Let $\{x_n\}$ be a Cauchy sequence in the Fr\'echet topology, and put
$x_n = \sum_i x_{n,i} t^i$. Then by hypothesis, for $i \geq 0$,
$\epsilon > 0$ and $r>1$,
there exists $N$ such that for $m, n \geq N$, $|x_{m,i} - x_{n,i}|r^i
< \epsilon$. In particular, for each fixed $i$, $\{x_{n,i}\}$ is a Cauchy
sequence in $\calO$, so it has a limit $y_i$.

Now put $y = \sum_i y_it^i$. For any given $\epsilon > 0$ and $r>0$,
choose $N$ so that for $m,n \geq N$, $w_r(x_m - x_n) < \epsilon$.
Then $|x_{m,i}-x_{n,i}| < \epsilon p^{-ri}$ for all $m,n \geq N$. Since
the absolute value on $\calO$ is nonarchimedean, we also have
$|x_{n,i} - y_i| < \epsilon p^{-ri}$, so $w_r(x_n - y) < \epsilon$. We
conclude that $x_n \to y$.
\end{proof}

We associate to each nonzero element $x$ of
$R$ a \emph{Newton polygon}
as follows. For each index $i$ with $x_i \neq 0$, plot the point
$(i, v_p(x_i))$ in the coordinate plane. We define the Newton polygon
of $x$ to be the
lower convex hull of these points and the \emph{slopes} of $x$ to be the
negations of the slopes of its Newton polygon (ignoring 0 if it occurs).
We define
the \emph{multiplicity} of a slope to be
the difference between the $x$-coordinates of the endpoints of the segment
of the Newton polygon with that slope (or 0 if there is no such segment).
We say $x \in R$ is \emph{pure} if $x$ has nonzero
constant coefficient
and exactly one slope; in particular, this implies $x \in \Oimm[\fp]$.

We single out a special class of pure elements that can be treated like
ordinary polynomials. We say $x = \sum x_i t^i$
is \emph{truncated} if $x$ is pure of some
slope $s$ with multiplicity $m$, and $x_i = 0$ for $i > m$.
%; observe that $w_{s'}(x) = p^{(s-s')m}$ for $s' < s$.
As for polynomials, there is a division lemma for truncated elements.
\begin{lemma}[Division lemma]
Let $x \in R$  be truncated of slope $s$ and multiplicity $m$.
Then for any $y \in R$, there exists a unique pair $q,r$ of elements
of $R$ such that $y = qx + r$ and $r = \sum r_i t^i$ satisfies $r_i=0$
for $i \geq m$. Moreover, $w_{s}(r) \leq w_{s}(y)$ and
$w_s(q) \leq w_s(y)$.
\end{lemma}
We will refer to $r$ as $y \bmod x$.
\begin{proof}
Since $x \in \Gimmcon[\fp]$, $x$ is invertible in that ring, and
$x^{-1}$ can be written as $t^{m} b$, where $b = \sum_i b_i t^i$ satisfies
$b_i=0$ for $i>0$. Now write $x^{-1}y = \sum_i z_i t^i$, and set
$q = \sum_{i \geq 0} z_it^i$. Then $x^{-1}y-q$ has no coefficients of
positive index, and so $r = y - qx = x(x^{-1}y-q)$ has no coefficients of
index $m$ or greater.

Note that $w_s(x^{-1})$ is well-defined; since $w_s(x) = 1$,
we must have $w_s(x^{-1}) = 1/w_s(x) = 1$. Thus $w_s(x^{-1}y) = 
w_s(y)$. Replacing some coefficients of a series with zeroes cannot
increase its norm, so $w_s(q) \leq w_s(x^{-1}y) = w_s(y)$
and $w_s(r) = w_s(x^{-1}r) = w_s(x^{-1}y - q) \leq w_s(x^{-1}y) 
\leq w_s(y)$.
\end{proof}

A \emph{slope factorization} of a nonzero element $x$ of $R$
is a product $x = t^j
\prod_{i=1}^{N} y_i$ for $N$ either a nonnegative
integer or $\infty$, convergent (if $N = \infty$) in the Fr\'echet topology,
such that:
\begin{enumerate}
\item[(a)]
$y_i$ is pure of slope $s_i$;
\item[(b)]
the sequence $s_1, s_2, \dots$ is strictly decreasing;
\item[(c)]
if $N = \infty$, then $s_i \to 0$.
\end{enumerate}
In our next few propositions,
we show that these factorizations exist and are unique up to units,
and use these factorizations to establish some structural properties
of $R$.

%PROPERTIES:
%Slope factors are unique up to units. BEGUN 
%Every element has a slope factorization. BEGUN
%Slope factorization is finite if and only if is in con subring.
%Every possible slope factorization occurs. BEGUN
%x divides y iff its slope factors do. BEGUN
%Newton polygon predicts $s_i$ (if we define the former!)

%We will use the Newton polygon of an element to give a canonical decomposition
%of the ideal generated by that element. Over $\Oan$, this decomposition
%can be expressed in terms of polynomials in $t$; for the larger rings $\Oalg$
%and $\Oimm$, we must define analogous entities and use those.

\begin{prop} \label{prop:newt}
Let $x$ and $y$ be elements of $R$ with nonzero constant coefficient. Then
the Newton polygon of $xy$ is the sum of the Newton polygons of $x$ and $y$.
That is, the multiplicity of a slope of $xy$ equals the sum of the
multiplicities of the corresponding slope of $x$ and of $y$.
\end{prop}
\begin{proof}
Fix a slope $s$. Suppose the Newton polygon of $x$ intersects its support
line of slope $-s$ from $(i, v_p(x_i))$ to $(j, v_p(x_j))$, and the Newton
polygon of $y$ intersects its support line of slope $-s$ from $(k, v_p(y_k))$
to $(l, v_p(y_l))$, with $i \leq j$ and $k \leq l$.
We claim the Newton polygon of $uv$ intersects its support line of slope
$-s$ at $(i+k, v_p(x_i)+v_p(y_k))$ and $(j+l, v_p(x_j) + v_p(y_l))$, which
would imply the statement of the lemma.

To verify the claim, we first note that
\begin{align*}
v_p((xy)_{i+k}) &=
v_p\left( \sum_m x_{i+m} y_{k-m} \right) \\
&\geq \min_{m \in \QQ} \{v_p(x_{i+m}) + v_p(y_{k+m})\},
\end{align*}
with equality if the minimum occurs for a single value of $m$.
For $m<0$, we have $v_p(x_{i+m}) > v_p(x_i) - ms$ and
$v_p(y_{k-m}) \geq v_p(y_k) + ms$, so 
$v_p(x_{i+m}) + v_p(y_k) > v_p(x_i) + v_p(y_k)$. Similarly,
if $m>0$, we have $v_p(x_{i+m}) \geq v_p(x_i) - ms$ and
$v_p(y_{k-m}) > v_p(y_k) + ms$, so again
$v_p(x_{i+m}) + v_p(y_k) > v_p(x_i) + v_p(y_k)$. Thus
the minimum is achieved only for $m=0$, and so
$v_p((xy)_{i+k}) = v_p(x_i) + v_p(y_k)$. Likewise,
$v_p((xy)_{j+l} = v_p(x_j) + v_p(y_l)$.

By a similar argument, we also have that for $n>0$,
$v_p((xy)_{i+k-n}) > v_p(x_i) + v_p(y_k)
+ ns$ and $v_p((xy)_{j+l+n} > v_p(x_j) + v_p(y_l) + ns$.
Namely,
\begin{align*}
v_p((xy)_{i+k-n}) &\geq \min_{m \in \QQ} \{ v_p(x_{i+m}) + v_p(y_{k-m-n}) \} \\
v_p(x_{i+m}) &\geq v_p(x_i) - ms \\
v_p(y_{k-m-n}) &\geq v_p(y_k) + (m+n) s,
\end{align*}
the second inequality is strict for $m<0$, and the third is strict
for $m>-n$, so $v_p(x_{i+m}) + v_p(y_{k-m-n}) > v_p(x_i) + v_p(y_k) + ns$
and the inequality follows; the other inequality follows analogously.
(The minimum really is a minimum, not an
infimum, so termwise strict inequality implies strict inequality for the
minimum.)
Thus the Newton polygon of $xy$ intersects its support line of slope $-s$
from $(i+k, v_p(x_i) + v_p(y_k))$ to $(j+l, v_p(x_j) + v_p(y_l))$,
respectively.
\end{proof}
\begin{cor}
The invertible elements of $\Gancon$ (resp.\ $\Galgancon, \Gimmancon$)
are precisely the nonzero elements
of $\Gcon[\fp]$ (resp.\ $\Galgcon[\fp], \Gimmcon[\fp]$).
\end{cor}

\begin{lemma} \label{lem:immfact0}
Let $x$ be an element of $R$ which has nonzero constant coefficient $1$.
Let $s$ be the first slope of $x$ and $m$ its multiplicity.
Then $x$ can be factored as $yz$ where $y$ is truncated of slope equal
to the first slope of $x$, with the same multiplicity.
\end{lemma}
\begin{proof}
We construct a sequence of elements $y_1, y_2, \dots$ of $R$ supported
on $[0, sm]$, having constant coefficient 1, convergent under the norm
$w_s$, such that $w_s(x \bmod y_n) \to 0$; this will imply that
$\{y_n\}$ and $\{x \bmod y_n\}$ converge in the Fr\'echet topology as well,
and that if $y_n \to 0$, then $x \bmod y = 0$. Specifically, we set
\begin{align*}
y_1 &= \sum_{0 \leq i \leq m} x_i t^i \\
y_{n+1} &= y_n - (x \bmod y_n) \qquad (n > 1).
\end{align*}
Let $c = w_r(x - y_1)$; by construction, $c < 1$. We prove by induction that
$w_r(x \bmod y_n) \leq c^n$. This holds by design for $n=1$. Now suppose
it holds up to some $n$. Then we can write $x = a_n y_n + b_n$ with
$b_n = x \bmod y_n$, and
\begin{align*}
x \bmod y_{n+1} &= (a_ny_n + b_n) \bmod y_{n+1} \\
&= (a_n (y_{n+1} - b_n) + b_n) \bmod y_{n+1} \\
&= b_n (1 - a_n) \bmod y_{n+1}.
\end{align*}
Since $w_r(y_n - y_1) \leq c$ by the induction hypothesis, and
$w_r(y_1) = 1$, we have $w_r(y_n) = 1$. Thus
\begin{align*}
w_r(1-a_n) &= w_r(y_n - a_n y_n) \\
&= w_r(y_n - x + b_n) \\
&\leq \max\{w_r(y_n - y_1), w_r(y_1 - x), w_r(b_n)\} \\
&\leq \max\{c, c, c^n\} = c.
\end{align*}
We conclude that $w_r(x \bmod y_{n+1}) \leq w_r(b_n(1-a_n)) \leq c^{n+1}$,
completing the induction.

By the induction, we have that $\{y_n\}$ is Cauchy, hence convergent,
and that $\{x \bmod y_n\}$ converges to 0. Thus the limit $y$ of $\{y_n\}$
satisfies $x \bmod y = 0$, as desired.
\end{proof}
\begin{cor}
Every pure element of $R$ factors as a truncated element times a unit.
\end{cor}

\begin{prop}
If $x$ has nonzero constant coefficient, then $(x,t)$ is the unit ideal.
\end{prop}
\begin{proof}
Let $x = \sum_i x_i t^i$. Without loss of generality, assume $x_0=1$,
and let $k$ be the smallest index such that
$v_p(x_k) < 0$. Let $y = x (\sum_{0 \leq i < k} x_i t^i)^{-1}$;
if we put $y = \sum_i y_i t^i$, then $y_i = 0$ for $0 < i < k$. That is,
$y \equiv 1 \pmod{t^k}$, so $(y, t^k)$ is the unit ideal, as then is
$(x, t)$.
\end{proof}

\begin{lemma} \label{lem:euclid}
\begin{enumerate}
\item[(a)] 
Let $x$ be a pure element of $R$ of slope $s$ and $y$ an element of $R$
with nonzero constant coefficient 
whose Newton polygon has all slopes greater than $s$. Then $(x,y)$
is the unit ideal.
\item[(b)]
Let $x$ and $y$ be pure elements of $R$ of slope $s$.
Then $(x,y)$ is generated by
a pure element of slope $s$.
\end{enumerate}
\end{lemma}
\begin{proof}
(a) Without loss of generality, we may assume $x$ is truncated and $y$
has constant coefficient 1.
Put $c = w_r(1-y)$; then $c<1$, and $w_r((1-y)^n) = c^n$.
Thus if we put $z_n = (1-y)^n \bmod x$, we also have $w_r(z_n) \leq c^n$,
so the series $\sum_{n=0}^\infty z_n$ converges in the Fr\'echet topology
to a limit $z$. Likewise, $\sum_{n=0}^\infty y(z_n \bmod x)$ converges
to 1, so $zy-1$ is divisible by $x$, and $(x,y)$ is thus the unit ideal.

(b) Let $m$ and $n$ be the multiplicities of $s$ as a slope of $x$ and $y$.
Then $m$ and $n$ are integral multiples of $v/s$, where $v$ is the smallest
positive valuation of $\calO$. Thus we can induct on $m+n$. Assume
without loss of generality that $m \leq n$, and that $x$ is truncated.
Let $z = y \bmod x$, so that $(x,y) = (x, z)$ and it suffices to show that
$(x,z)$ is generated by a pure element of slope $s$.

If $z=0$, then $(x,y) = x$ and we are done, so assume $z \neq 0$.
Otherwise, let $z = \sum_i z_i t^i$ and choose $j$ to minimize
$v_p(z_j) + rj$. Since $(x,t)$ is the unit ideal by the previous
lemma, we can find $a$ such that $az_jt^j \equiv 1 \pmod{x}$, and 
$(x,z) = (x, az)$. Let $b = az \bmod x = \sum_i b_i t^i$; then $b_0 = 1$ and
all slopes of $b$ are at least $s$. Moreover, the multiplicity of
$s$ as a slope of $b$ is strictly less than $m$.
By Lemma~\ref{lem:immfact0}, we can factor $b$ as $cd$ with $c$ pure of
slope $s$ with the same multiplicity, and $d$ having all slopes greater
than $s$. By (a), $(x,d)$ is the unit ideal, so $(x,z)$, which is
equal to $(x,cd)$, is also equal to $(x,c)$. By the induction hypothesis,
$(x,c)$ is generated by a pure element of slope $s$, as desired.
\end{proof}

\begin{prop} \label{prop:divis}
For $x,y \in R$ such that $x$ admits a slope factorization,
$x$ divides $y$ if and only if each factor in a slope
factorization of $x$ divides $y$.
\end{prop}
\begin{proof}
If $x$ divides $y$, then obviously any factor of $x$ divides $y$.
Conversely, suppose $ct^j \prod y_i$ is a slope factorization of $x$.
Put $z_i = y / (ct^j y_1\cdots y_i)$; then the sequence $\{z_i\}$ is Cauchy,
so has a limit $z$, and clearly $zx \to y$.
\end{proof}

\begin{prop} \label{prop:conv}
  Let $\{z_n\}_{n=1}^\infty$ be a sequence of pure elements of $R$ whose
slopes are strictly increasing and tend to $0$. Then there exists $x \in R$
admitting a slope factorization $c t^j \prod y_n$ such that $y_n$ and $z_n$
generate the same ideal for all $n$.
\end{prop}
\begin{proof}
Let $s_n$ be the slope of $z_n$, and let $v$ be the smallest positive
valuation of $\calO$. Put $z_n = \sum_i z_{n,i} t^i$; then $z_{n,i} \in 
\calO$ for $i < v/s_n$. In particular, $\sum_{i < v/s_n} z_{n,i} t^i$
is a unit in $R$; let $u_n$ be its inverse and put $y_n = u_n z_n$.
If we put $y_n = \sum_i y_{n,i} t^i$, then $y_{n,i} = 0$ for $0<i<v/s_n$
by construction and $v_p(y_{n,i}) \geq -s_n i$ for $i \geq v/s_n$.
In particular, we have $w_r(y_n-1) \leq p^{-v-rv/s_n}$ for
$r \geq s_n$. Therefore for any fixed $r>1$, eventually $w_r(y_n-1)
\leq p^{-v-rv/s_n}$; since $s_n \to 0$, $p^{-rv/s_n} \to 0$.
We conclude that $\prod_n y_n$ converges in the Fr\'echet topology,
and we may take the limit as our desired $x$.
\end{proof}

\begin{lemma} \label{lem:immfact2}
Every nonzero element $x$ of $R$
has a slope factorization.
\end{lemma}
\begin{proof}
Without loss of generality, assume $x$ has nonzero constant coefficient.
Factor $x$ as $y_1 x_2$ as in Lemma~\ref{lem:immfact0}, so that $y_1$
is pure of slope equal to the first slope of $x$, with the same multiplicities.
Then by Proposition~\ref{prop:newt}, the Newton polygon of $x_2$ is equal
to that of $x$ with its first segment removed. Then factor $x_2$ as
$y_2 x_3$ in the same fashion, and so on.

If the Newton polygon of $x$ has finitely many slopes, then this process
eventually represents $x$ as a product of pure elements, as desired.
If the Newton polygon of $x$ is infinite, Proposition~\ref{prop:conv}
allows us to construct an element $z$ with a slope factorization
whose slope factors are unit multiples of the $y_i$. In particular,
$x$ and $z$ have the same Newton polygon. Since $x$ is divisible by
each $y_i$, by Proposition~\ref{prop:divis} $x$ is divisible by $z$.
Moreover, the Newton polygon of $x/z$ is empty, so $x/z$ is a unit.
We can modify the slope factorization of $z$ by multiplying its first factor
by $x/z$ to obtain the desired slope factorization of $x$.
\end{proof}

The next proposition
is a rigid analytic version of the Chinese remainder theorem.
\begin{prop}[Chinese remainder theorem] \label{prop:crt}
Let $\{x_n\}_{n=1}^\infty$ be a sequence of pure elements of $R$ 
whose slopes are strictly increasing and
tend to $0$, and let $\{y_n\}_{n=1}^\infty$ be elements of $R$.
Then there exists $y \in R$ such that $y \equiv y_n \pmod{x_n}$.
\end{prop}
\begin{proof}
Let $s_n$ be the slope of $x_n$, and let $v$ be the smallest positive
valuation in $\calO$. By imitating the proof of Proposition~\ref{prop:conv},
we may assume that $x_n = \sum_i x_{n,i} t^i$ is such that $x_0=1$ and
$x_{n,i} = 0$ for $0<i<v/s_n$. In particular, this implies that
$\prod x_n$ converges to a limit $x$; let $u_n = x/x_n$.

We construct a sequence $\{z_n\}_{n=1}^\infty$ such that
$u_n z_n \equiv y_n \pmod{x_n}$ and $\sum u_n z_n$ converges;
then we may set $y = \sum u_n z_n$ and be done.
First, choose $v_n$ with $u_n v_n \equiv y_n \pmod{x_n}$, which exists
because $u_n$ and $x_n$ are relatively prime.

It suffices to show that
for any $\epsilon>0$, there exists $c$ such that for all $r \geq p^{s_{n-1}}$,
$w_r(1+cx_n)
< \epsilon$; for example, one can then choose such an $\epsilon$ with
$\epsilon w_r(u_nv_n) < 1/n$ for all $r \geq s_{n-1}$,
and set $z_n = v_n(1+cx_n)$.

Note that $w_r(1-x_n) < 1$ for $r > p^{s_n}$, so in those norms,
the sequence $c_m = -1 -(1-x_n) - \cdots - (1-x_n)^m$ is Cauchy. In
particular, there exists $m$ such that $w_r(u_n v_n(1 + cx_n)) < \epsilon$
for all $r \geq s_{n-1}$. Set $z_n = v_n(1+cx_n)$; it is now clear
that the series $\sum u_n z_n$ is Cauchy for all of the $w_r$, and thus
convergent. We then take $y= \sum u_nz_n$ and the proof is complete.
\end{proof}

\begin{prop}
The slope factorization of $x$ is unique up to units.
\end{prop}
\begin{proof}
The constant coefficient and power of $t$ are clearly uniquely determined
by $x$, so we may suppose $\prod y_i$ and $\prod z_j$ are slope factorizations
of $x$. It suffices to prove that each $z_j$ divides one of the $y_i$ (and
vice versa). Taking the greatest common divisor of $z_j$ with the $y_i$
of the same slope, if it exists, and dividing out that divisor, we may
reduce to the case where each $y_i$ is relatively prime to $z_j$.

In this case, for each $i$, there exists $w_i$ such that $w_iz_j \equiv 1
\pmod{y_i}$. By Proposition~\ref{prop:crt}, there exists $w$ such that
$wz_j \equiv 1 \pmod{y_i}$ for each $i$, so by Proposition~\ref{prop:divis}
$wz_j-1$ is divisible by $x$. But then $wz_j$ and $wz_j-1$ are both divisible
by $z_j$, a contradiction since $z_j$ is not a unit.
Thus $z_j$ divides one of the $y_i$, as desired.
\end{proof}

Finally, we use slope factorization to prove a structure theorem about ideals
in analytic rings.
\begin{lemma} \label{lem:principal}
Every finitely generated ideal in $\Gancon$ (resp.\ $\Galgancon, \Gimmancon$) is principal.
\end{lemma}
\begin{proof}
It suffices to show that for $x, y \in \Gancon$ (resp.\ $x,y \in \Gimmancon$),
there exist $r,s$ such that $rx+sy$ generates the ideal $(x,y)$.
By Lemma~\ref{lem:immfact1}, we may reduce to the case where $x,y \in \Oan$
(resp.\ $\Oalg_{\an}, \Oimm_{\an}$) and have nonzero constant coefficients.
By Proposition~\ref{prop:conv}, there exists $z$ whose slope factorization
consists of the greatest common divisors of the slope factors of $x$ and $y$,
and by Proposition~\ref{prop:divis}, $x$ and $y$ are divisible by $z$. Dividing
off $z$, we reduce to the case where $x,y$ have no common slope factors.

Let $x = c \prod y_i$ be the slope factorization of $x$. By assumption,
$y$ is not divisible by $y_i$, so there exists $z_i$ such that $yz_i \equiv 1
\pmod{y_i}$. By Proposition~\ref{prop:crt}, there exists $z$ such that
$z \equiv z_i \pmod{y_i}$; then $yz-1$ is divisible by all of the $y_i$.
Therefore there exists $w$ such that $yz-1 = wx$, and $1 \in (x,y)$ as desired.
\end{proof}

\section{A direct calculation}

The critical non-formal step in the proof of the main theorem
is a direct computation to establish the existence of certain
eigenvectors. This computation takes the form of two related lemmas.

\begin{lemma} \label{lem:special}
Let $M$ be an $F$-crystal over $R = \Galgancon$ or $R = \Gimmancon$,
and let $\lambda \in \calO_0$ be a uniformizer.
Suppose $M$ admits a basis $\bv_1, \dots, \bv_n, \bw$
such that for some
$c_i \in R$,
\begin{align*}
F\bv_i &= \bv_{i+1} \qquad (i=1, \dots, n-1) \\
F\bv_n &= \lambda^{n+1} \bv_1 \\
F\bw &= \bw + \sum_{i=1}^n c_i \bv_i.
\end{align*}
Then $M$ has an eigenvector $\by$ with $F\by = \lambda \by$.
\end{lemma}
\begin{proof}
Observe that for $j=1, \dots, n-1$,
\[
F(\bw + c \bv_j) = (\bw+c\bv_j) + c \bv_{j+1} - \bv_j + \sum_{i=1}^n c_i \bv_i.  
\]
Thus we can modify $\bw$ by a suitable linear combination of $\bv_1,
\dots, \bv_n$ to obtain $\bx$ such that $F\bx = \bx + y \bv_1$
for some $y \in R$.

Suppose $\by = a \bx + b_1 \bv_1 + \cdots + b_n \bv_n + b \bx$ satisfies
$F\by = \lambda \by$.
Comparing coefficients in this equation, we have
$a^\sigma = \lambda a$,
$b_i^{\sigma} = \lambda b_{i+1}$ for $i=1, \dots, n-1$,
and $\lambda^{n+1} b_n^\sigma + ax = b_1$.
If $a$ and $b$ satisfy the equations
\begin{equation} \label{eq:system}
a^\sigma = \lambda a, \qquad
\lambda b^{\sigma^n} = b - \lambda^{-1} a x,
\end{equation}
then setting $b_1 = b$ and $b_i = b_1^{\sigma^{i-1}} \rho^{-i+1}$ for
$i=2, \dots, n$ and $\by = a\bx + b_1 \bv_1 + \cdots + b_n \bv_n
+ b\bx$ satisfies $F\by = \lambda \by$.
%\fixme{check arithmetic}
Thus it suffices to show that (\ref{eq:system}) has
a solution.

Notice that replacing $x$ by $x + \lambda^{n+1} y^{\sigma^n} - y$
does not alter whether a solution exists: for any $a$
such that $a^\sigma = \lambda a$,
\[
a (\lambda^{n+1} y^{\sigma^n} - y) = 
\lambda^{-n} [ \lambda^{n+1} (ay)^{\sigma^n} - \lambda^n ay].
\]
If we write $x = c + d$ with $c = \sum_{i<0} x_i t^i$
and $d = \sum_{i\geq 0} x_i t^i$, then the equation
$\lambda^{n+1} y^{\sigma^{n}}
- y = d$ has a solution by Lemma~\ref{lem:rankone}.
Thus without loss of generality, we may assume $x_i = 0$ for $i \geq 0$. In
particular, we now have $x \in \Gimmcon[\fp]$.

We next reduce to the case where $x$ is supported on $[-1, -1/p^{n})$.
Set
\[
y = \sum_{j=1}^\infty \sum_{k=1}^j \sum_{i \in [-1, -1/p^m)}
\lambda^{(n+1)(k-1)} x_{ip^{-nj}}^{\sigma^{m(k-1)}}
+ \sum_{j=1}^\infty \sum_{k=1}^j \sum_{i \in [-1, -1/p^m)}
\lambda^{-(n+1)k} x_{ip^{nj}}^{\sigma^{-mk}}.
\]
The first sum is evidently convergent in $R$. As for the second,
recall that there exists constants $c,d$ such that $v_p(x_{-n})
\geq cn-d$ for $n < 0$. Thus
\[
v_p(\lambda^{-(n+1)k} x_{ip^j}) \geq (-i)cp^j - d
- j(n+1) v_p(\lambda),
\]
and the right side grows exponentially in $j$, so the second
sum is also convergent.
Thus we can replace $x$ by
\[
x - \lambda^{n+1} y^{\sigma^{n}}
+ y
= \sum_{i \in [-1, -1/p^{n}]} \sum_{j=-\infty}^\infty 
x_{ip^{-nj}}^{\sigma^{nj}} \lambda^{(n+1)j},
\] 
which is supported on $[-1, -1/p^{n})$.

We now assume $x$ is supported on $[-1, -1/p^n)$.
If $x=0$, then $M$ has an eigenvector $\bx$ with $F\bx = \bx$, which
can be multiplied by a suitable scalar to produce the desired $\by$,
so we assume $x \neq 0$; indeed, we may assume $|x|=1$.
Define $c_i$ for all $i < 0$ by setting $c_i = x_i$ for
$i \in [-1, -1/p^m)$ and extending by the rule $c_{ip^n} = \lambda^{n+1}
c_i^{\sigma^n}$. We say an index $i$ is a \emph{corner} of $c$
if $|c_j| < |c_i|$ for all $j<i$; then all corners are of the form
$jp^{-nk}$ for $j$ in a finite set and $k$ an arbitrary integer.

The solutions of $a^{\sigma} = \lambda a$ have the form
\[
a = \sum_{i \in [1, p)}
\sum_{j=-\infty}^\infty a_{i}^{\sigma^{j}}
\lambda^{-j} t^{ip^{j}}
\]
for $\sum_{i \in [1, p)} a_i t^i \in R$.
For such $a$, define
\begin{align*}
f(a) &= \sum_{i\in [l, p^{-n}l)} t^i
\sum_{j=-\infty}^\infty (xa)_{ip^{-nj}}^{\sigma^{nj}} \lambda^{j} \\
&= \sum_{i \in [l, p^{-n}l)} t^i
\sum_{j=-\infty}^\infty \sum_{k \in [-1,-1/p^n)}
x_k^{\sigma^{nj}} a_{ip^{-nj}-k}^{\sigma^{nj}} \lambda^{j} \\
&= \sum_{i \in [l, p^{-n}l)} t^i
\sum_{j=-\infty}^\infty \sum_{k \in [-1,-1/p^n)}
c_{kp^{nj}} a_{i-kp^{nj}} \\
&= \sum_{i \in [l, p^{-n}l)} t^i
\sum_{k<0} c_k a_{i-k}.
\end{align*}
Then $f$ is additive, and
$f(a) = 0$ if and only if the equation $\lambda b^{\sigma^n}
= b - \lambda^{-1} a x$ has a solution.

For $i > 0$ and $\alpha \in \calO$,
put $a(\alpha,i) = \sum_{j=-\infty}^\infty
\alpha^{\sigma^j} \lambda^{-j} t^{ip^j}$. 
Write 
\[
f(a(\alpha, i)) = \sum_{s \in [l, p^{-n}l)}
t^s \sum_{j=\infty}^\infty \lambda^{-j} c_{s-ip^j}
\alpha^{\sigma^j}.
\]
Let $j(i)$ be the smallest $j$ which achieves
$\max_{s,j} \{|\lambda^{-j} c_{s-ip^j}|\}$
for some $s$, let $s(i)$ be the smallest such $s$,
and let $k(i) = s(i) - ip^{j(i)}$.
Then $k(i)$ is a corner of $c$ if $s(i) \neq l$.
On the other hand, if $s(i) = l$ and 
$k(i)$ were not a corner, we could find $k' < k$ with
$|c_{k'}| \geq |c_k|$. Let $m>0$ be the unique integer such that
$p^{-mn}(k' - ip^{j(i)}) \in [l, pl)$. Then
\[
|\lambda^{-j(i)+mn} c_{p^{-mn} k' - ip^{j(i)-mn}}|
= |\lambda^{-j(i)-m} c_{k'-ip^{j(i)}}|
> |\lambda^{-j(i)} c_{k'-ip^{j(i)}}|,
\]
contradiction. Thus $k(i)$ is also a corner if
$s(i) = l$, so $k(i)$ is piecewise constant, as is $j(i)$;
therefore $s(i)$ is piecewise linear and increasing in $i$.
Also, clearly $s(ip) = s(i)$ and $j(ip) = j(i)-1$.

We claim that if $\lim_{i' \to i^-} s(i')$ and
$\lim_{i' \to i^+} s(i')$ lie in $(l, p^{-n}l)$
for some $i$, then $s$ is continuous at $i$.
Let $k_0$ and $k_1$
be the value of $k(i-\epsilon)$ and $k(i+\epsilon)$, respectively,
for small $\epsilon>0$; define $j_0$ and $j_1$ analogously.
If $|\lambda^{-j_0} c_{k_0}| < |\lambda^{-j_1} c_{k_1}|$, then
$|\lambda^{-j_0} c_{s-(i+\epsilon)p^{j_0}}| < |\lambda^{-j_1} c_{k_1}|$
for $s = k_0 + (i+\epsilon)p^{j_0}$, contradicting the definition
of $j(i+\epsilon)$. A similar contradiction arises if
$|\lambda^{-j_0} c_{k_0}| < |\lambda^{-j_1} c_{k_1}|$. Hence
$\lambda^{-j_0} c_{k_0}$ and $\lambda^{-j_1} c_{k_1}$ have the same norm.
Now by the definition of $s$, we have $k_0 +(i-\epsilon)p^{j_0}
\leq k_1 + (i-\epsilon)p^{j_1}$ and
$k_0 + (i+\epsilon)p^{j_0} \geq k_1 + (i+\epsilon)p^{j_1}$.
Thus $k_0 + ip^{j_0} = k_1 + ip^{j_1}$ and $s$ is continuous at $i$.
%moreover, $|c_{k_0}/c_{k_1}| = |\lambda^{j_0-j_1}|$.

We can also show that $\lim_{i' \to i^+} s(i') = l$ if and only if
$\lim_{i' \to i^-} s(i') = p^{-n}l$, but the argument is more delicate.
Define $k_0, k_1, j_0, j_1$ as above. 
If $\lim_{i' \to i^-} s(i') = p^{-n} l$ but
$\lim_{i' \to i^+} s(i') = s_1 > l$,
 then
for $i' = i - \epsilon$ we have
$|\lambda^{-j} c_{k_1}| =
|\lambda^{-j} c_{s_1-ip^{j_1}}| < |\lambda^{-j_0} c_{k_0}|$.
On the other hand, for $i' = i+\epsilon$, we can take
$s = p^{n}(k_0 + i'p^{j_0}) < s_1$ and obtain
\[
|\lambda| |\lambda^{-j_0} c_{k_0}|
=
|\lambda^{-j_0-n} c_{k_0p^{n}}|
< |\lambda^{-j_1} c_{k_1}|.
\]
Thus $|\lambda^{-j_1} c_{k_1}|$ is sandwiched between
$|\lambda||\lambda^{-j_0}c_{k_0}|$ and $|\lambda^{-j_0}c_{k_0}|$,
but there is no norm between these two because $\lambda$ is a uniformizer,
contradiction. If $\lim_{i' \to i^+} s(i') = l$ but
$\lim_{i' \to i^-} s(i) = s_0 < p^{-n}l$,
then $|\lambda^{-j_0} c_{k_0}| = |\lambda^{-j_0} c_{s_0-ip^{j_0}}|
\leq |c_{k_1} \lambda^{-j_1}|$.
On the other hand, for $i' = i - \epsilon$,
we can take $s = p^{-n}(k_1 + (i-\epsilon)p^{j_1}) > s_0$ and obtain
\[
|\lambda^{-1}| |\lambda^{-j_1} c_{k_1}|
=
|\lambda^{-j_1+n} c_{k_1p^{-n}}|
\leq |\lambda^{-j_0} c_{k_0}|
\]
but $|\lambda^{-1}| > 1$, so again we obtain a contradiction. 

Since $s(ip) = s(i)$ and $s$ is increasing and
continuous except for jumps between
$p^{-n}l$ and $l$, it follows that $s$ maps $[1, p)$ onto $[l, p^{-n}l)$
one or more times. Choose $r \in [1, p)$ such that
$s(r) = l$; then $s$ also maps $[r, pr)$ onto $[l, p^{-n}l)$ one or more
times. It follows that
for any $y$ supported on $[l, p^{-n}l)$
with $|y| \leq 1$, there exists $a$ with 
$a^\sigma = \lambda a$ with $|a_i| \leq 1$ for $i \in [r,pr)$
such that $|y-f(a)| < 1$.
This $a$ can be constructed by a transfinite recursion: find $i \in [r,pr)$
and $\alpha \in R$ such that
$f(a(\alpha,i))$ has leading coefficient congruent to the leading
coefficient of $y$ modulo $\lambda$, then subtract
off and repeat.

Additionally, note that $s$ must change slope
at some point in $[r, pr)$ (possibly equal to $r$), since
$j(r) \neq j(pr)$. If $s$ changes slope at $i$, then
$f(a(\alpha, i))$ has at least two distinct terms with minimal norm.
Namely, if again $k_0$ and $j_0$ (resp.\ $k_1$ and $j_1$) are the
values of $k(i-\epsilon)$ and $j(i-\epsilon)$ (resp.\ $k(i+\epsilon)$
and $j(i+\epsilon)$) for $\epsilon>0$ small, then the terms
$c_{k_0} \lambda^{-j_0} \alpha^{\sigma^{j_0}}$ and
$c_{k_1} \lambda^{-j_1} \alpha^{\sigma^{j_1}}$ have the same minimal
norm. Note that there exists $\alpha \in \calO$ such that the sum of the
terms of minimal norm has norm less than 1 (because the residue field
of $\calO$ is algebraically closed). Thus in the transfinite recursion
of the previous paragraph, there is more than one choice that can be made
at $i$. In particular,
there exists $a$ with $a^\sigma = \lambda a$,
$|f(a)|<1$ and $|a_i| = 1$ for some $i \in [r, pr)$.

From this analysis, we can construct a nonzero solution of (\ref{eq:system}).
Start with $a^{(0)}$ such that $(a^{(0)})^\sigma = \lambda a^{(0)}$,
$|f(a^{(0)})| < 1$ and $|a^{(0)}_1| = 1$ for some $i \in [r, pr)$.
Now construct a sequence $\{a^{(m)}\}_{m=0}^\infty$ such that 
\begin{enumerate}
\item[(a)] $(a^{(m)})^\sigma = \lambda a^{(m)}$ for all $m$;
\item[(b)] $f(a^{(m)}) < |\lambda^m|$;
\item[(c)] $|a^{(m)}_i - a^{(m+1)}_i| < \lambda^m$ for $i \in [r, pr)$.
\end{enumerate}
Specifically, given $a^{(m)}$, let $y = f(a^{(m)})/\lambda^m$,
find $a$ such that $a^\sigma = \lambda a$ with $a_i \leq 1$ for $i \in [r,pr)$
and $|y-f(a)| < 1$, then set $a^{(m+1)} = a^{(m)} + \lambda^m a$.
Then $\{a^{(m)}\}$ converges in the Fr\'echet topology to a nonzero
$a$ with $a^\sigma = \lambda a$ and $f(a) = 0$. Thus
(\ref{eq:system}) has a nonzero solution, and the proof is complete.

\end{proof}

\begin{lemma}
  \label{lem:twod}
Let $M$ be an $F$-crystal over $R = \Galgancon$ or $R = \Gimmancon$
admitting a basis $\bv, \bw$ such that
\begin{align*}
  F\bv &= \pi^m \bv \\
F\bw &= \bw + c \bv
\end{align*}
for $\pi\in \calO_0$ a uniformizer, $c \in R$,
and $m \geq 2$. Then there exists an eigenvector $\bx$ of $M$ such tha
$F\bx = \pi^{m-1} \bx$.
\end{lemma}
\begin{proof}
We imitate the previous proof restricted to $n=1$.
For starters, given $x \in R$ supported on $[-1, -1/p)$ with $|x|=1$,
it again suffices to show that
the equations
\[
a^\sigma = \pi^{m-1} a, \qquad
\pi b^{\sigma} = b - \pi^{-1} a x
\]
have a solution with $a,b \in R$ not both zero (the case $x=0$ is
self-evident).
Define $c_i$ for all $i < 0$ by setting $c_i = x_i$ for
$i \in [-1, -1/p^m)$ and extending by the rule $c_{ip^n} = \lambda^{n+1}
c_i^{\sigma^n}$. Again, we say an index $i$ is a \emph{corner} of $c$
if $|c_j| < |c_i|$ for all $j<i$.

For $a$ a solution of $a^\sigma = \pi^{m-1} a$,
define
\[
f(a) = \sum_{i \in [l, l/p)} t^i
\sum_{k<0} c_k a_{i-k};
\]
it suffices to exhibit $a$ such that $f(a) = 0$.

Continuing to imitate the previous proof,
for $i>0$, define $j(i)$, $k(i)$, $s(i)$ 
by taking $j(i)$ as the smallest $j$ for which
$\max_{s,j} \{|\pi^{-(m-1)j} c_{s-ip^j}|\}$ is achieved,
$s(i)$ as the smallest $s$ for which the maximum is achieved
with $j=j(i)$, and $k(i) = s(i) - ip^{j(i)}$.
Then as before, $k(i)$ and $j(i)$ are piecewise constant and so
$s(i)$ is increasing and piecewise linear.

We again prove that $\lim_{i' \to i^+} s(i')$ and $\lim_{i' \to i^-}
s(i')$ are either equal, or equal to $l$ and $l/p$, respectively.
The proof that if both limits lie in $(l, l/p)$, then they are equal,
carries through as before,
as does the proof that we cannot have $\lim_{i' \to i^+} s(i') = l$
and $\lim_{i' \to i^-} s(i') < l/p$.
Now suppose $\lim_{i' \to i^-} s(i') = l/p$ but
$\lim_{i' \to i^+} s(i') = s_1 > l$.
For $i' = i - \epsilon$ we have
$|\pi^{-(m-1)j} c_{k_1}| =
|\pi^{-(m-1)j} c_{s_1-ip^{j_1}}| < |\pi^{-(m-1)j_0} c_{k_0}|$.
On the other hand, for $i' = i+\epsilon$, we can take
$s = p(k_0 + i'p^{j_0}) < s_1$ and obtain
\[
|\pi| |\pi^{-(m-1)j_0} c_{k_0}|
=
|\pi^{-(m-1)(j_0+1)} c_{k_0p}|
< |\pi^{-(m-1)j_1} c_{k_1}|.
\]
Thus $|\pi^{-(m-1)j_1} c_{k_1}|$ is sandwiched between
$|\pi||\pi^{-(m-1)j_0}c_{k_0}|$ and $|\pi^{-(m-1)j_0}c_{k_0}|$,
but there is no norm between these two because $\pi$ is a uniformizer,
contradiction.

Given the results of the previous paragraph, the rest of the proof
proceeds as in the previous lemma. Namely, a transfinite recursion can
be used to generate a solution of $f(a)=0$, which completes the proof.

\end{proof}

\section{Construction of the special Newton polygon}

This section is devoted to the proof of the following theorem, the main result
of this paper. For $n$ a positive integer and $s$ a rational number such that
$sn$ is a valuation of an element $\lambda \in \calO$,
let $M_{n,s}$ denote the crystal over $\Omega$
whose action of Frobenius on a basis $\bv_1, \dots, \bv_n$ is given by
$F\bv_i = \bv_{i+1}$ for $i=1, \dots, n-1$ and $F\bv_n = \lambda \bv_1$.
We say a crystal is \emph{standard} if it is isomorphic to some $M_{n,s}$.
\begin{theorem} \label{thm:main}
Let $R = \Galgancon$ or $R = \Gimmancon$.
Let $M$ be an $F$-crystal over $R$ spanned by eigenvectors.
Then $M$ splits as a direct sum of standard subcrystals.
\end{theorem}
\begin{cor} \label{cor:special}
Let $M$ be an $F$-crystal over $R = \Galgcon$ or $R = \Gimmcon$.
Then $M$ becomes constant when extended to $R_{\an} \otimes_{\calO}
\calO'$ for some suitable finite extension $\calO'$ of $\calO$.
\end{cor}

By Lemma~\ref{lem:slopes} below, the slopes of the eigenvectors
that form a basis of $M$ do not depend on the choice of basis. We call these
the \emph{special slopes} of $M$,
and we define the \emph{special Newton polygon} of $M$ as
the convex polygon whose slopes are the special slopes of $M$; we
will catalog its basic properties in the next section.

%Beware that crystals over $\Gimmancon$ have more eigenvectors
%than one might expect; in fact, if a crystal has a basis of eigenvectors,
%then it has eigenvectors of every sufficiently large slope (over a
%suitable extension of $\calO$). However, the following lemmas shows that
%the additional assumption that an eigenvector is part of a basis for
%the crystal cuts down the number of candidates significantly.

\begin{lemma} \label{lem:slopes}
Suppose $D$ and $E$ are diagonal matrices over $\calO$ and $U$ is an
invertible matrix over $\Gimmancon$ such that $U^{-1} D U^{\sigma} = E$.
Then the slopes (valuations of the diagonal entries) of $D$ and $E$ are the
same up to a permutation.
\end{lemma}
\begin{proof}
Since $U$ has nonzero determinant, we can find a permutation matrix $V$
such that $UV$ has nonzero diagonal entries. Put $W = UV$ and 
$F = V^{-1}EV$; then $F$ is diagonal and its entries are a permutation
of those of $E$.
From $DV^{\sigma} = VF$ we have $D_{ii} V_{ii}^\sigma = V_{ii} F_{ii}$
for each $i$; since $V_{ii}$ is nonzero, this implies $|D_{ii}| \leq
|F_{ii}|$. In particular, the $k$-th largest slope of $D$ is greater than
or equal to the $k$-th largest slope of $E$ for each $k$. The analogous
statement with $D$ and $E$ reversed follows by a similar argument; therefore
the slopes of $D$ and $E$ are equal up to permutation.
\end{proof}

For $\bv$ an element of an $F$-crystal over a ring $R$, the ideal generated
by the coordinates of $\bv$ in some basis is independent of the choice of
basis; we call this ideal the \emph{coordinate ideal} of $\bv$. We
say $\bv$ is \emph{primitive} if its coordinate ideal is the trivial ideal.
Equivalently, $\bv$ is primitive if and only if it can be extended to a 
basis of $R$.
\begin{lemma} \label{lem:primitive}
Let $M$ be an $F$-crystal over $R_{\an}$ for $R = \Galgcon$ or $R=\Gimmcon$
and $\bv \in M$ a nonzero
element such that $F\bv = \mu\bv$. Then $\bv$ is a multiple of a primitive
eigenvector of $M$. In particular, if $M$ has no eigenvectors
of slope less than $v_p(\mu)$, then $\bv$ is primitive.
\end{lemma}
\begin{proof}
Let $I$ be the coordinate ideal of $\bv$.
By Lemma~\ref{lem:principal}, $I$ is principal,
so we may choose a generator $r$. Since $F\bv = \mu \bv$, the ideal $I$
is invariant under $\sigma$ and $\sigma^{-1}$,
so $r^\sigma = cr$ for $c$ a unit in $R$. 
By Lemma~\ref{lem:immfact1} we can write $c = \lambda d$ with
$\lambda \in \calO$ and $d$ a unit in $R$.
The equation $t^\sigma = dt$ has a solution with $t \in R$,
and $s = r/t$ satisfies $s^\sigma = \lambda s$ and generates
$I$. Therefore there exists $\bw \in M$ with $\bv = s \bw$,
$\bw$ primitive, and $F\bw = (\mu/\lambda)\bw$, as desired.
\end{proof}

Beware that there may be primitive eigenvectors which cannot be extended
to a basis consisting solely of eigenvectors. For example, there is always
an eigenvector of slope equal to the largest generic slope, which usually
does not extend to a basis of eigenvectors. However, it will turn out that
the eigenvector of minimum slope will always extend to a basis of eigenvectors.

\begin{proof}[Proof of Theorem~\ref{thm:main}]
We proceed by induction on $n = \dim M$.
For $n=1$, $F$ acts on a basis vector by an invertible scalar, i.e.\
an element of $\Gimmcon[\fp]$. Without loss of generality, we may assume
this scalar is in $\Gimmcon$ and has norm 1, in which case the result
follows from \cite[Lemma~2.1]{bib:me4}.

Now suppose $n>1$.
We are done if $M$ is isomorphic to $M_{n,d/n}$, so we assume that this
does not occur.
Let $d$ be the slope of $\wedge^n M$.
For $r$ a rational number, define the \emph{$\calO$-index} of $r$ as the
smallest positive integer $s$ such that $rs$ is a valuation of $\calO$.

Since the
set of rationals of $\calO$-index less than $n$ is discrete,
there is a smallest such rational that occurs as the slope of an
eigenvector of $M$ over a suitable extension of $\calO$;
call this number $r$.
Let $m$ be the $\calO$-index of $r$,
and let $\lambda \in \calO_0$
have valuation $rm$. Let $\bv$ be an eigenvector of $M$ over $R[\lambda^{1/m}]$
with $F\bv = \lambda^{1/m} \bv$. Write $\bv = \sum_{i=0}^{m-1}
\bw_i \lambda^{-i/m}$, so that each $\bw_i$ is an element of $M$ over $R$
with $F^m \bw_i = \lambda \bw_i$, and let $M_1$ be the span of $\bw_0,\dots,
\bw_{m-1}$ within $M$. 

Since $\dim M_1 \leq m < n$, we may
apply the induction hypothesis to $M_1$. If $M_1$ has more than one
standard summand, it has an eigenvector of slope less than or equal to $r$
with $\calO$-index strictly less than $m$. Thus $M$ has an eigenvector
of slope strictly less than $r$ with $\calO$-index less than $n$,
contradiction.
Thus $M_1$ itself is standard; specifically, it must
be isomorphic to $M_{m,r}$ (in particular, $\dim M_1$ must equal $m$).
Likewise, $M/M_1$ has a direct sum decomposition $N_2 \oplus \cdots \oplus
N_k$ of the specified type.

Let $P_i$ be the preimage of $N_i$ in $M$; to complete the proof,
it suffices to show that the exact sequence
$0 \to M_1 \to P_i \to N_i \to 0$ splits for $i=2, \dots, k$.
First suppose $k>2$. Then the dimension of each $P_i$ is less than $n$,
so we can apply the induction hypothesis to $P_i$.
If the slope of $N_i$ were less than that of $M_1$, then 
the induction hypothesis would imply that $P_i$ has a slope 
less than $r$ of $\calO$-index less than or equal to $\dim P_i \leq n$,
contradicting the minimality of $r$. Thus
the slope of $N_i$ is greater than or equal to
that of $M_1$. In this case,
Lemma~\ref{lem:rankone} can be used to show that the exact
sequence splits. Namely, by imitating the argument at the beginning
of the proof of Lemma~\ref{lem:special}, we can
reduce this splitting to the
existence of solutions of equations of the form $\lambda a^{\sigma^d} - 
\mu a = x$, where $d$ is the least common multiple of
$\dim M_1$ and $\dim N_1$, and $|\lambda| > |\mu|$. Then
Lemma~\ref{lem:rankone} implies that each of these equations
has a solution.

The case $k=2$ requires special scrutiny, as $P_i = M$ and the
induction hypothesis does not apply.
Let $s$ be the slope of $N_2$.
As above, the exact sequence splits if $r \leq s$, so assume
on the contrary that $r > s$; this implies in particular that $r > d/n$.
If $n=2$, we immediately obtain a contradiction from Lemma~\ref{lem:twod},
so we may assume $n>2$.

We show that
$M$ has an eigenvector of slope less than or equal to
$d/n$ over some finite extension of $\calO$. 
Pick an eigenvector $\bv$ of slope $r$, and apply the induction hypothesis
to the quotient of $M$ by the span of this eigenvector. This yields
an eigenvector $\bw$ of the quotient of slope at most $(d-r)/(n-1)$. Applying
the induction hypothesis again, this time to the preimage of $\bw$,
gives an eigenvector of $M$ of some slope $r' \leq (r + (d-r)/(n-1))/2$.
Let $\calO^1$ be a finite extension of $\calO$ whose value group
contains $r'$, and let $r_1$ be the smallest rational of
$\calO^1$-index less than $n$ that occurs as the slope of an eigenvector.
Again, take the quotient of $M$ by the span of an eigenvector of slope
$r_1$, this time over $\calO^1$ and apply the induction hypothesis. If
its slopes are not all equal, we deduce that $M$ has the desired splitting
over $\calO^1$,
and in particular has an eigenvector of slope less than or equal to
$d/n$ over $\calO^1$. Otherwise, we can repeat the process to
produce an extension $\calO^2$ of $\calO$, and the smallest rational
$r_2$ of $\calO^2$-index less than $n$ occuring as the slope of an
eigenvector will be at most $(r_1 +(d-r_1)/(n-1))/2$, and so on.

The existence of an eigenvector of slope at most $d/n$ is assured if
the above process ever terminates, so assume it continues indefinitely.
Then the sequence $\{r_n\}$ converges to $d/n$, as it is sandwiched
between $d/n$ and a sequence obtained by iterating
$r \mapsto (r + (d-r))/(n-1))/2$, and the latter converges to $d/n$;
therefore, there exist eigenvectors of $M$ of every rational slope
greater than $d/n$.
Let $\calO'$ be an extension of $\calO$ whose value group contains
$d/n$, and let $v_0$ be its minimum positive valuation. Take
an eigenvector $\bv$ of $M$ of slope $d/n + v_0/(n-1)$ over an
extension of $\calO'$ of degree $n-1$. The span of $\bv$ over 
$\calO$ has dimension at most $n-1$ and sum of slopes at most
$(n-1)d/n + v_0$. Apply the induction hypothesis to the span;
if the sum of slopes is not equal to $(n-1)d/n+v_0$, then
it is at most $(n-1)d/n$, and
the span contains an eigenvector of slope at most $d/n$. If
the sum of slopes is equal to $(n-1)d/n+v_0$ but the slopes are not
all equal, and $t$ is the least slope and its multiplicity is $m$,
then $mt < md/n+mv_0/(n-1)$; since $mt$ and $md/n$ are integral
multiples of $v_0$, we deduce $mt \leq md/n$ and the span again contains
an eigenvector of slope at most $d/n$. Finally, if all of the slopes
of the span are equal to $d/n+v_0/(n-1)$, then
Lemma~\ref{lem:special} implies that $M$ has an eigenvector of slope
$d/n$. Thus in all cases, $M$ has an eigenvector of slope less than 
or equal to $d/n$.

Let $\lambda$ be an element of a finite
extension of $\calO_0$ of valuation $d/n$. Because $M$ has an eigenvector
of slope less than or equal to $d/n$, it must also have one of slope
equal to $d/n$, over some finite extension of $\calO[\lambda]$.
In fact, from a solution of $F\bv = \lambda \bv$ over a finite extension
of $\calO[\lambda]$, we can obtain a solution over $\calO[\lambda]$: decompose
the solution over a basis of the finite extension over $\calO[\lambda]$,
and choose any nonzero component. Thus we may assume $\bv$ is 
defined over $\calO[\lambda]$.

Let $m$ be the $\calO$-index of $d/n$; then we can write
$\bv = \sum_{i=0}^{m-1} \lambda^{-i} \bw_i$. Let $N$ be the
span of $\bw_0, \dots, \bw_{m-1}$. If $\dim N < n$, then we can apply
the induction hypothesis to $N$ to express it as a direct sum of
standard subcrystals. The projection of $\bw_0$ onto at least one of these
subcrystals must be nonzero, and since $F^m \bw_0 = \lambda^m \bw_0$,
the same equation holds for the projection of $\bw_0$. Thus this subcrystal
has an eigenvector of slope $d/n$; its slope must then be
less than $d/n$. Since this slope has $\calO$-index at most $m \leq n$,
we conclude $r \leq d/n$, contradiction. On the other hand, if $\dim N = n$,
then $N$ is isomorphic to $M_{n,d/n}$. Moreover,
$\bw_0 \wedge \cdots \wedge \bw_{n-1}$ is an eigenvector of 
$\wedge^n M$ of slope $d$, so must be primitive. Thus $N=M$ is isomorphic
to $M_{n,d/n}$, contrary to an earlier assumption. We conclude that
the assumption $r > s$ leads to a contradiction in all cases, so we must
have $r \leq s$ and so the exact sequence
$0 \to M_1 \to P_2 \to N_2 \to 0$ splits.

In summary, we have that each $P_i$ splits as a direct sum of $M_1$ with
another summand; now $M$ splits as a direct sum of these other summands
with $M_1$, as desired.
\end{proof}

We do not know whether or not an arbitrary $F$-crystal over $\Gimmancon$
is spanned by eigenvectors, and hence has a basis of eigenvectors.
Indeed, if it were known that
every $F$-crystal over $\Gimmancon$ has a nonzero eigenvector, it would
follow by induction
that every $F$-crystal is spanned by eigenvectors.

\section{Properties of the Newton polygons}

Unless otherwise specified throughout
this section, let $M$ be an $F$-crystal over $\Gimmcon$.
In this section, we establish that the special and generic Newton polygons
satisfy some relations that one would expect from the case of potentially
semistable $(F, \nabla)$-crystals over $\Gcon$.

\begin{prop}
Let $M$ be an $F$-crystal over $\Omega$. Then the special Newton polygon
of $M$ equals the Newton polygon
of the reduction of $M$ modulo $t$.
\end{prop}
This follows immediately from Dwork's trick \cite[Lemma~4.3]{bib:me3}.
Beware that the natural generalization of this
proposition to $\Oimm$ is false.

\begin{prop}
If $\ell_1 \leq \cdots \leq \ell_n$ are the special slopes of $M$,
the special slopes of $\wedge^k M$ (for $k = 0, \dots, m$) are
given by $\ell_{i_1}+\cdots+\ell_{i_k}$ for $1 \leq i_1 < \cdots < i_k \leq m$.
\end{prop}
\begin{proof}
If $\bv_1, \dots, \bv_n$ form a basis of eigenvectors of $M \otimes_{\Gimmcon}
\Gimmancon$ with $F\bv_i = \lambda_i \bv_i$, then
\[
F(\bv_{i_1}\wedge\cdots \wedge \bv_{i_k}) =
\lambda_{i_1} \cdots \lambda_{i_k} \bv_{i_1} \wedge \cdots \bv_{i_k},
\]
so $\bv_{i_1} \wedge \cdots \wedge \bv_{i_k}$ is an eigenvector of 
slope $v_p(\lambda_{i_1}) + \cdots + v_p(\lambda_{i_k})$.
\end{proof}
Similarly, the special slopes of $M_1 \oplus M_2$ are the union of
the special slopes of $M_1$ and $M_2$, and the special slopes of
$M_1 \otimes M_2$ are the products of the special slopes of $M_1$
and $M_2$.

\begin{prop}
The special Newton polygon lies above the generic Newton 
polygon and has the same endpoints.
\end{prop}
\begin{proof}
By the previous proposition, it suffices to prove that the highest
special slope is no greater than the highest generic slope.
Let $\bv_1, \dots, \bv_n$ be a basis of eigenvectors of $M \otimes_{\Gimmcon}
\Gimmancon$ and let $\bw$ be an eigenvector of $M$ of highest slope over
$\Gimmcon$, which exists by \cite[Proposition~2.1]{bib:me4}.
Write $\bw = \sum c_i \bv_i$ with $c_i \in \Gimmancon$.
If $F\bw = \lambda \bw$ and $F\bv_i = \mu_i \bv_i$, then for $i=1, \dots, n$,
we have
$c_i^\sigma \mu_i = c \lambda$, so $c_i^\sigma = (\lambda/\mu_i) c_i$. The
latter equation only has solutions in $\Gimmancon$ if $|\lambda/\mu_i| \leq 1$.
Thus each special slope is less than or equal to the highest generic slope.
\end{proof}

\begin{prop}
The vertices of the special Newton polygon are of the form $(i, j)$, where
$i$ is an integer and $j$ is an integral multiple of the smallest positive
valuation in $\calO$.
\end{prop}
\begin{proof}
This follows immediately from Theorem~\ref{thm:main}.
\end{proof}

%Properties:
%Special over generic
%Stable under base extension
%Correct for semistable crystals

\begin{prop}
Let $M$ be an $F$-crystal over $\Gcon$ whose generic and special
Newton polygons coincide. Then $M$ becomes unipotent over $\Gsepcon
\otimes_{\calO} \calO'$ for an extension $\calO'$ of $\calO$ whose
value group contains the slopes of $M$.  
\end{prop}
\begin{proof}
It suffices to show that the eigenvectors of $M$ of lowest slope are 
defined over $\Gsepcon$.
Let $\bv_1, \dots, \bv_n$ be eigenvectors
of $M$ over $\Galgancon$, with $F\bv_i = \lambda_i \bv_i$ for
$\lambda_i \in \calO_0$ such that $|\lambda_1| \leq \cdots \leq |\lambda_n|$.
By the descending slope filtration for $F$-crystals
\cite[Lemma~2.1]{bib:me4},
we can find elements $\bw_1, \dots, \bw_n$ of $M$ over $\Galgcon$ such that
$F\bw_i = \lambda_i \bw_i + \sum_{j<i} A_{ij} \bw_j$ for some
$A_{ij} \in \Galgcon$.

Write $\bv_n = \sum_i b_i \bw_i$ with $b_i \in \Galgcon$;
then applying $F$ to both sides,
we have $\lambda_n b_i = \lambda_i b_i^\sigma + \sum_{j>i} b_j^\sigma A_{ji}$
for $i=1, \dots, n$. For $i=n$, this implies $b_n \in \calO_0$.
By Lemma~\ref{lem:rankone} and descending induction on $i$,
we obtain $b_i \in \Galgcon$ for $i=n-1, \dots, 1$, and so
$\bv_n$ is defined over $\Galgcon$. 
On the other hand, by the ascending slope filtration
\cite[Proposition~2.2]{bib:me4},
the eigenvectors $\bv$ of $M$ over $\Galg$ satisfying
$F\bv = \lambda_n \bv$ are defined over $\Gsep$. Thus 
$\bv_n$ is defined over $\Galgcon \cap \Gsep = \Gsepcon$.
Since $\bv_n$ could have been taken to be any eigenvector over
$\Galgancon$ of lowest slope, this proves the claim.
\end{proof}
\begin{cor} \label{cor:filtr}
Let $M$ be an $F$-crystal over $\Gcon$ whose generic and special
Newton polygons coincide. Then $M$ has a filtration
$0 = M_0 \subseteq M_1 \subseteq \cdots \subseteq M_l = M$ by subcrystals
such that $M_i/M_{i-1}$ is isoclinic for $i=1, \dots, n$.
\end{cor}

%This statement can be refined for $(F, \nabla)$-crystals using the
%fact that \'etale $(F,\nabla)$-crystals become constant over
%a finite separable extension of Tsuzuki \cite{bib:tsu1}.
%The result is a special case of the semistable
%reduction conjecture \cite[Conjecture~4.12]{bib:me3}.
%\begin{cor}
%Let $M$ be a $(F, \nabla)$-crystal over $\Gcon$
%whose generic and special Newton polygons coincide. Then $M$ becomes
%unipotent over a finite separable extension of $\Gcon$.
%\end{cor}
%\begin{proof}
%It suffices to show that the eigenvectors of lowest slope, which
%are defined over $\Gsepcon$, are in fact defined over a finite
%separable extension of $\Gcon$; applying this assertion to the exterior
%powers of $M$ proves the claim.
%The coincidence of the Newton polygons is equivalent to $M$ becoming
%a constant $F$-crystal over $\Galgcon$.
%In fact, the eigenvectors of lowest slope are annihilated by $\nabla$,
%since an eigenvector of slope $s$ is sent to an eigenvector of slope $s-1$.
%Thus the lowest slope eigenspace of
%$M$ is a sub-$(F,\nabla)$-crystal over $\Gcon$;
%by Tsuzuki's theorem \cite[Theorem~5.1.1]{bib:tsu1}, its eigenvectors
%are defined over a finite extension of $\Gcon$, as desired.
%\end{proof}

\section*{Acknowledgments}
This work is partially based on the author's doctoral dissertation
\cite{bib:methesis}, written under the supervision of Johan de~Jong.
The author was supported by a National Science Foundation Postdoctoral
Fellowship. Thanks also to Laurent Berger for helpful discussions.

\end{document}